\documentclass[10pt,reqno]{amsart}
  \usepackage{geometry}
\newtheorem*{theorem A}{Theorem A}
\newtheorem*{theorem B}{N\"olker's Theorem}

\theoremstyle{remark}

\theoremstyle{remark}

\theoremstyle{definition}

\numberwithin{equation}{section}
\def\({\left ( }
\def\){\right )}
\def\<{\left < }
\def\>{\right >}

\begin{document}

\title[ON THE EXISTENCE OF PARA-KENMOTSU MANIFOLDS]
 {ON THE EXISTENCE OF\\ PARA-KENMOTSU MANIFOLDS}

\author[BELDJILALI Gherici]{BELDJILALI Gherici}

\address{Laboratory of Quantum Physics and Mathematical Modeling (LPQ3M),\\ University of Mascara,  Algeria}
\email{gherici.beldjilali@univ-mascara.dz}



\date{June 06, 2018}


\begin{abstract}
This note provides a quite obvious observation that the condition (2.7), which is a part of the original definition of the so-called para-Kenmotsu manifolds \cite{SS1}, does not make sense, and thus this concept is void. So, it is proved that the para-Kenmotsu manifolds does not exist under the condition mentioned above.
\end{abstract}

\maketitle
\section{Introduction}
The notion of almost para-contact structure was introduced by Sato \cite{IS1, IS2}. This structure is an analogue of the almost
contact structure \cite{BL1, YK} and is closely related to almost product structure (in contrast to almost contact structure, which is related to almost complex structure). An almost contact manifold is always odd dimensional but an almost para-contact manifold could be even dimensional as well.

After that, T. Adati and K. Matsumoto \cite{AM} defined and studied para-Sasakian and SP-Sasakian manifolds which are regarded as a special kind of an almost contact Riemannian manifolds. Before Sato, Kenmotsu \cite{KEN} defined a class of almost contact Riemannian manifolds.

Later, in 1995, Sinha and Sai Prasad \cite{SS1} have defined a class of almost para contact metric manifolds namely para-Kenmotsu (p-Kenmotsu) and special para Kenmotsu (sp-Kenmotsu) manifolds. Unfortunately without example which to prove the existence of this type of manifolds. 

Para-Kenmotsu manifolds in this sense, are recently developed   in several works, for example as  \cite{SS2, SS3, SSD}.

This note provides a quite obvious observation that the formula (2.7), which is a part of the original definition of the so-called para-Kenmotsu manifolds \cite{SS1}, does not make sense, and thus this concept is void. 
\\
Firstly, let's give a brief information for para-Kenmotsu manifolds.
\section{Preliminaries}
Let $M$ be an $n$-dimensional differentiable manifold equipped with structure tensors $( \varphi , \xi, \eta)$ where $\varphi$ is a
tensor of type $(1, 1)$, $\xi$ is a vector field, $\eta$ is a 1-form such that:
\begin{equation}\label{eq1}
\eta(\xi) = 1,
\end{equation}
\begin{equation}\label{eq2}
\varphi^2 X = X - \eta(X)\xi,
\end{equation}
for all $X$ vector field on $M$. Then $M$ is called an almost para-contact manifold.\\
Let $g$ be the Riemannian metric in an $n$-dimensional almost para-contact manifold $M$ such that
\begin{equation}
\eta(X) = g(\xi , X),
\end{equation}
\begin{equation}
\varphi \xi =0, \qquad \eta(\varphi X)=0,\qquad rank \varphi = n-1,
\end{equation}
\begin{equation}\label{eq3}
g(\varphi X , \varphi Y) =g(X,Y)-\eta(X)\eta(Y),
\end{equation}
for all vector fields $X$ and $Y$ on $M$. Then the manifold $M$  is said to admit an almost para-contact Riemannian
structure $( \varphi , \xi, \eta , g)$ and the manifold is called an almost para-contact Riemannian manifold.\\

A manifold $M$ with Riemannian metric $g$ admitting a tensor field $\varphi$ of type $(1, 1)$, a vector field $\xi$ and 1-form $\eta$
satisfying equations (\ref{eq1}), (\ref{eq2}), (\ref{eq3}) along with
\begin{equation}\label{eq4}
{\rm d}\eta =0, \qquad \nabla_X \xi = \varphi^2 X = X- \eta(X)\xi,
\end{equation}
\begin{eqnarray}\label{eq5}
\big( \nabla_X \nabla_Y \eta \big) Z =& \big( -g(X,Z)+\eta(X)\eta(Z) \big)\eta(Y)\notag \\
& +   \big( -g(X,Y)+\eta(X)\eta(Y) \big)\eta(Z),
\end{eqnarray}
is called a para-Kenmotsu manifold or briefly P-Kenmotsu manifold \cite{SS1}, where $\nabla$ denote the Levi-Civita connection
with respect to $g$.

In the following, we will  give three reasons why the condition (\ref{eq5}) is never satisfied.

\section{Existence of para-Kenmotsu manifolds}

Our objection to the condition (\ref{eq5}) lies in the following three serious problems:\\

\textbf{First problem:}\\

In the basic condition (\ref{eq5}), putting
$$ T(X,Y)Z=\big( -g(X,Z)+\eta(X)\eta(Z) \big)\eta(Y) +   \big( -g(X,Y)+\eta(X)\eta(Y) \big)\eta(Z).$$
Then,  for all $f$ function on $M$, we have

\begin{eqnarray*}
 \big( \nabla_X \nabla_{fY} \eta \big) Z& =& \Big( \nabla_X \big( f \nabla_Y\eta \big)\Big) Z \\
 &=&  X (f) \big( \nabla_Y\eta \big) Z +  f \big( \nabla_X  \nabla_Y\eta \big) Z.
 \end{eqnarray*}
But in the other side we have
 \begin{eqnarray*}
T(X,fY)Z&=&\big( -g(X,Z)+\eta(X)\eta(Z) \big)\eta(fY) +   \big( -g(X,fY)+\eta(X)\eta(fY) \big)\eta(Z) \\
 &=&  f \Big(\big( -g(X,Z)+\eta(X)\eta(Z) \big)\eta(Y) +   \big( -g(X,Y)+\eta(X)\eta(Y) \big)\eta(Z)\Big)\\
 &=& f T(X,Y)Z,
 \end{eqnarray*}
i.e. 
 \begin{equation*}
 \big( \nabla_X \nabla_{fY} \eta \big) Z = f T(X,Y)Z + X (f) \big( \nabla_Y\eta \big) Z. \\
 \end{equation*}
 
This means that one side of the condition (\ref{eq5}) is $C^{\infty}(M)$-linear in the second variable $Y$ the other  is not.\\
 
 \textbf{Second problem:}\\

In the basic condition (\ref{eq5}),  since ${\rm d}\eta =0$ i.e for all $Y$ and $Z$ vector fields on $M$ 
$$ (\nabla_Y \eta ) Z = (\nabla_Z \eta ) Y.$$
 Then, we have

\begin{eqnarray*}
 \big( \nabla_X \nabla_Y \eta \big) Z& =&\nabla_X  \Big( \nabla_Y\eta \Big) Z  -  \nabla_Y\eta \Big( \nabla_X  Z \Big)\\
 &=& \nabla_X  \Big( \nabla_Z\eta \Big) Y  -  \nabla_Y\eta \Big( \nabla_X  Z \Big)\\
 &\neq &  \big( \nabla_X \nabla_Z \eta \big) Y,
 \end{eqnarray*}
and in the other side we have
 \begin{eqnarray*}
T(X,Y)Z &=&\big( -g(X,Z)+\eta(X)\eta(Z) \big)\eta(Y) +   \big( -g(X,Y)+\eta(X)\eta(Y) \big)\eta(Z) \\
 &=& T(X,Z)Y,
 \end{eqnarray*}
This means that one side of the condition (\ref{eq5})  admits the commutativity between $Y$ and $Z$ the other  is not.\\
 
 \textbf{Third problem}:\\
 
 From (\ref{eq4}), we have
  \begin{eqnarray*}
 2{\rm d}\eta(X,Y) =0 &\Leftrightarrow& \big( \nabla_X \eta \big) Y =  \big( \nabla_Y \eta \big) X\\
  &\Leftrightarrow& g(\nabla_X \xi , Y) = g( \nabla_Y \xi , X),
  \end{eqnarray*}
  by setting $Y = \xi$, we find
   \begin{equation}\label{eq6}
\nabla_{\xi} \xi =0.
 \end{equation}
Therefore, from  (\ref{eq5}), we have

 \begin{equation}\label{eq7}
\big( \nabla_X \nabla_{\xi} \eta \big) Z =-g(X,Z)+\eta(X)\eta(Z) ,
\end{equation}
and  in the other side we have
 \begin{eqnarray}\label{eq8}
 \big( \nabla_X \nabla_{\xi} \eta \big) Z& =&  \nabla_X \Big( \big( \nabla_{\xi} \eta \big)Z\Big) -  \big( \nabla_{\xi} \eta \big) \big(\nabla_X Z \big) \notag\\
  &=&  X g\big( \nabla_{\xi} \xi  , Z\big) -  g\big( \nabla_{\xi} \xi  , \nabla_X Z\big) \notag\\
 &=&  0.
 \end{eqnarray}
 From (\ref{eq7}) and (\ref{eq8}) shows the contradiction.\\

 This clearly shows that condition  (\ref{eq5}) is never fulfills and therefore the definition of para-Kenmotsu structure as well as the theorem 3.1 in \cite{SDH} must be reconsidered.

\end{document}